\documentclass[12pt]{amsart}

\usepackage{enumitem}
\usepackage{psfrag}
\usepackage{graphicx}
\usepackage{pinlabel}
\usepackage{graphicx}
\usepackage{amsmath}
\usepackage{amssymb}
\usepackage{amscd}
\usepackage{times}
\usepackage{pb-diagram}
\usepackage{graphicx}
\usepackage{wrapfig}
\usepackage{xspace}
\usepackage{hyperref}
\usepackage{url}
\usepackage{caption}
\usepackage{subcaption}
\usepackage{fancyhdr}
\usepackage{overpic}
\usepackage{color}
\usepackage[square,comma,sort&compress,numbers]{natbib}
\usepackage{latexsym}
\usepackage{mathrsfs}
\usepackage{amsfonts}
\usepackage{hyperref}
\usepackage{amsthm}
\usepackage{appendix}
\usepackage{pdfsync}

\theoremstyle{definition}
\newtheorem{theorem}{Theorem}[section]
\newtheorem{lemma}[theorem]{Lemma}
\newtheorem{proposition}[theorem]{Proposition}

\newtheorem{definition}[theorem]{Definition}

\newtheorem{remark}[theorem]{Remark}
\newtheorem*{theorem*}{Theorem}
\def\qed{\hfill{Q.E.D.}\smallskip}

\title{Combinatorial Calabi Flows on Infinitely Triangulated Surfaces}
\author{Yuerong Bian, Xu Xu}
\date{\today}

\address{School of Mathematics and Statistics, Wuhan University, Wuhan, 430072, P. R. China}
\email{yuerongbian@whu.edu.cn}

\address{School of Mathematics and Statistics, Wuhan University, Wuhan, 430072, P. R. China}
\email{xuxu2@whu.edu.cn}

\thanks{MSC (2020): 52C25, 52C26}

\keywords{Vertex scaling; combinatorial Calabi flow; infinitely triangulated surfaces; piecewise Euclidean metrics; piecewise hyperbolic metrics}

\begin{document}
	
	\begin{abstract}
		We study the combinatorial Calabi flow associated with vertex scalings
		on infinitely triangulated surfaces with piecewise Euclidean and
		piecewise hyperbolic metrics.
		Under the assumptions of uniformly bounded vertex degree and an
		$\epsilon$-uniformly nondegenerate initial metric, we establish
		short-time existence of smooth solutions and prove uniqueness among
		solutions that are continuous in the $L^\infty$ norm.

	\end{abstract}

\maketitle

\section{Introduction}

\subsection{Background}

Geometric flows provide a method for studying metrics with prescribed
curvature.
Hamilton introduced the Ricci flow \cite{Hamilton}, while Calabi
introduced the Calabi flow to study extremal K\"ahler metrics
\cite{CALABI4,CALABI5}.
The long-time existence and convergence of the Calabi flow on closed
Riemannian surfaces have been studied in
\cite{CHANG6,CHANG7,CHEN8,CHRUSCIEL11}.

On triangulated surfaces, combinatorial curvature flows are used to
study polyhedral metrics with prescribed combinatorial curvature.
Chow and Luo \cite{BL} introduced the combinatorial Ricci flow for
Thurston's circle packing metrics \cite{THURSTON29} on closed triangulated surfaces  
and studied its long-time existence and convergence.
Luo \cite{Luo1} introduced the combinatorial Yamabe flow for vertex
scaling of piecewise Euclidean metrics. 
Ge \cite{Gphd,GE13} introduced the combinatorial Calabi flow for Thurston's circle packing metrics and vertex scalings on surfaces.
Related developments include combinatorial Yamabe flows with surgery for vertex scalings
\cite{GLSW,GGLSW}, combinatorial Calabi flows with surgery for vertex scalings \cite{ZXu}, 
combinatorial Ricci flows with surgery for inversive distance circle packings \cite{BL1, BL2},
combinatorial Calabi flows with surgery for inversive distance circle packings \cite{XZ3, XZ4},
fractional combinatorial Calabi flows \cite{WuXu} and others.

The preceding results on combinatorial curvature flows concern
finitely triangulated closed surfaces.
For infinitely triangulated surfaces, Ge, Hua and Zhou \cite{GHZ}
studied combinatorial Ricci flows for Thurston's circle packing
metrics on infinite disk triangulations in Euclidean and hyperbolic
background geometries.
Ji \cite{JI} investigated the combinatorial Yamabe flow for Euclidean
vertex scalings on infinitely triangulated surfaces.
In \cite{BX}, we studied the hyperbolic combinatorial Yamabe flow
on infinitely triangulated surfaces and established short-time
existence and uniqueness under the assumptions of uniformly bounded
vertex degree and an $\epsilon$-uniformly nondegenerate initial metric.

In this paper, we study the combinatorial Calabi flow for vertex
scalings on infinitely triangulated surfaces in both Euclidean
and hyperbolic background geometries.
Under the assumptions of uniformly bounded vertex degree and
an $\epsilon$-uniformly nondegenerate initial metric, we establish
short-time existence of smooth solutions, with existence times
depending only on $\epsilon$ and the degree bound.
We also prove uniqueness among solutions that are continuous
in the $L^\infty$ norm; the solutions constructed in the
existence proof satisfy this continuity condition.

The combinatorial Calabi flow requires additional estimates beyond those
used for the hyperbolic combinatorial Yamabe flow in \cite{BX}.
For the combinatorial Yamabe flow, a uniform curvature bound directly
controls the time derivative of the conformal factor.
The right-hand side of the combinatorial Calabi flow  also contains the
derivatives of the curvature with respect to the conformal
factors.
We establish uniform bounds for these derivatives and their
variation in both background geometries.
These estimates control the first and second time derivatives
of the approximating solutions, with bounds independent of
the finite subcomplex.
Using an exhaustion argument as in Ge, Hua and Zhou
\cite{GHZ}, we obtain solutions on a common time interval
during which the metrics remain uniformly nondegenerate.
For uniqueness, we use uniform Lipschitz estimates and
Gr\"onwall's inequality.
This argument has been applied to the Ollivier Ricci
flow on infinite graphs by Hua, Lin and Liu
\cite[proof of Theorem~1.2]{HLL}, and to the hyperbolic
combinatorial Yamabe flow in our previous work \cite{BX}.
For the combinatorial Calabi flow, we additionally need uniform estimates
for the differences of the curvature derivatives with
respect to the conformal factors.
Together with the curvature estimates, these yield the
required Lipschitz estimate for the right-hand side
of the flow in the $L^\infty$ norm.

\subsection{Definitions and main results}

Let $S$ be a connected surface without boundary, and let
$\mathcal T=(V,E,F)$ be a locally finite triangulation of $S$.
Throughout this paper, the vertex set $V$ is assumed to be infinite.
We use $i$, $\{ij\}$, and $\triangle ijk$ to denote a vertex, an edge,
and a face, respectively.
Two vertices $i,j$ are adjacent, denoted by $i\sim j$, if
$\{ij\}\in E$.
The number of edges incident to $i$ is called its degree and is
denoted by $\deg(i)$.

A piecewise Euclidean metric, also called a piecewise linear
(PL) metric, is a flat cone metric on $S$ whose cone singularities
are contained in $V$.
A piecewise hyperbolic (PH) metric is a hyperbolic cone metric on
$S$ whose cone singularities are contained in $V$.
In either case, the metric is determined by a length function
$d:E\to(0,\infty)$ such that the edge lengths of each face satisfy
the strict triangle inequalities.
These lengths determine a nondegenerate Euclidean or hyperbolic
triangle, according to the background geometry.
The metric is obtained by gluing the corresponding triangles
isometrically along their common edges.
We identify a PL or PH metric with its edge length function $d$.

The combinatorial curvature at a vertex $i\in V$ is defined by
\[
K_i=2\pi-\sum_{\triangle ijk\in F}\theta_i^{jk},
\]
where the sum is taken over the triangles incident to $i$, and
$\theta_i^{jk}$ is the interior angle at $i$ in the corresponding
Euclidean or hyperbolic triangle.

Vertex scaling gives a notion of discrete conformality for
polyhedral metrics on triangulated surfaces.
For PL metrics, it was introduced independently by
Ro\v{c}ek--Williams \cite{Rocek} and Luo \cite{Luo1}.
Bobenko, Pinkall and Springborn \cite{BPS} introduced the
corresponding notion for PH metrics.
We recall the definitions below.

\begin{definition}\label{Def: vertex_scaling}
	Let $d$ and $\widetilde d$ be two metrics on $(S,\mathcal T)$.
	\begin{enumerate}
		\item[(a)] If $d$ and $\widetilde d$ are PL metrics, 
		$\widetilde d$ is called a Euclidean vertex scaling of $d$ if there is a function $u:V\to\mathbb R$ such that
		\[
		\widetilde d_{ij}=e^{(u_i+u_j)/2}d_{ij},
		\qquad \forall \{ij\}\in E.
		\]
		\item[(b)] If $d$ and $\widetilde d$ are PH metrics, 
		$\widetilde d$ is called a hyperbolic vertex scaling of $d$ if there is a function $u:V\to\mathbb R$ such that
		\[
		\sinh\frac{\widetilde d_{ij}}2
		=e^{(u_i+u_j)/2}\sinh\frac{d_{ij}}2,
		\qquad \forall \{ij\}\in E.
		\]
	\end{enumerate}
	In either case, $u$ is called a discrete conformal factor, and we
	write $\widetilde d=u\ast d$.
\end{definition}

Vertex scaling on closed surfaces has been extensively studied.
See, for example, Bobenko-Pinkall-Springborn \cite{BPS} on the global rigidity of Euclidean and hyperbolic vertex scalings, 
Gu-Luo-Sun-Wu \cite{GLSW} and Gu-Guo-Luo-Sun-Wu
\cite{GGLSW} for discrete uniformization theorems for PL and
PH metrics, respectively.

In this paper, 
we consider the following combinatorial Calabi flows on a surface $S$ with fixed infinite
triangulation $\mathcal T$.

\begin{definition}\label{Def:VCLBF}
	Let $d_0$ be an initial PL or PH metric on $(S,\mathcal T)$.
	The combinatorial Calabi flow associated with vertex scalings is
	defined by
	\begin{equation}\label{Eq:VCLBF}
		\begin{cases}
			\dfrac{du_i}{dt}=\Delta K_i, & i\in V,\\[6pt]
			u_i(0)=0, & i\in V,
		\end{cases}
	\end{equation}
	where the discrete Laplace operator is given by
	\begin{equation}
		\Delta f_i
		=-\sum_{j\in V}\frac{\partial K_i}{\partial u_j}f_j
	\end{equation}
	for any function $f:V\to\mathbb R$.
	Here $K_i(u)$ is the combinatorial curvature of $u\ast d_0$ at $i$,
	and all partial derivatives are evaluated at $u$.
\end{definition}

The flow is considered as long as $u(t)\ast d_0$
remains a PL or PH metric of the same type as $d_0$.
Since $\mathcal T$ is locally finite and
$\partial K_i/\partial u_j=0$ whenever $j\ne i$ and $j\not\sim i$,
the sum defining $\Delta f_i$ contains only finitely many
nonzero terms for each fixed vertex $i$.

We use the following $\epsilon$-uniform nondegeneracy
condition on the initial metric.

\begin{definition}
	Let $d$ be a PL or PH metric on $\mathcal T$.
	For $\epsilon>0$, we say that $d$ is
	$\epsilon$-uniformly nondegenerate if
	\[
	\min\{\theta_i^{jk},\theta_j^{ik},\theta_k^{ij}\}
	\geq\epsilon
	\]
	for every triangle $\triangle ijk\in F$.
\end{definition}

If $0<\epsilon_1\leq\epsilon_2$, then every
$\epsilon_2$-uniformly nondegenerate metric is also
$\epsilon_1$-uniformly nondegenerate.

We denote by $L^\infty(V)$ the space of bounded functions on $V$,
with norm
\[
\|u\|_{L^\infty}=\sup_{i\in V}|u_i|.
\]
The notation $u\in C_t^\infty(V\times[0,T])$ means that
$u_i\in C^\infty([0,T])$ for every $i\in V$.
Continuity in $C([0,T];L^\infty(V))$ is understood with respect to the above norm.

In this paper, we prove the following short-time existence result for the combinatorial Calabi flow (\ref{Eq:VCLBF}) on infinitely triangulated surfaces.
\begin{theorem}[Short-time existence]\label{thm:vxste}
	Let $(S,\mathcal T)$ be an infinitely triangulated surface with
	vertex degree bounded by $D$, and let $d_0$ be an
	$\epsilon$-uniformly nondegenerate initial metric on $\mathcal T$.
	\begin{enumerate}
		\item[(a)] If $d_0$ is a PL metric, there exists
		$T_1=T_1(\epsilon,D)>0$ such that the flow \eqref{Eq:VCLBF}
		admits a solution
		$
		u\in C_t^\infty(V\times[0,T_1]).
	$
		\item[(b)] If $d_0$ is a PH metric, there exists
		$T_2=T_2(\epsilon,D)>0$ such that the flow \eqref{Eq:VCLBF}
		admits a solution
	$
		u\in C_t^\infty(V\times[0,T_2]).
	$
	\end{enumerate}
\end{theorem}

The solutions constructed in the proof of
Theorem~\ref{thm:vxste} are also continuous in time with respect
to the $L^\infty$ norm. This follows from the uniform bounds
for $|du_i/dt|$ established in Sections~2 and~3.
Under the same assumptions on the infinitely triangulated surface
and the initial metric, we further prove the following short-time uniqueness
result for solutions with this continuity property.

\begin{theorem}[Short-time uniqueness]
	\label{thm:VCLBF_uniqueness}
	Let $(S,\mathcal T)$ be an infinitely triangulated surface
	with vertex degree bounded by $D$, and let $d_0$ be an
	$\epsilon$-uniformly nondegenerate initial PL or PH metric.
	Let $T>0$ satisfy
	\[
	T\leq
	\begin{cases}
		T_1(\epsilon,D), & \text{if $d_0$ is a PL metric},\\
		T_2(\epsilon,D), & \text{if $d_0$ is a PH metric},
	\end{cases}
	\]
	where $T_1$ and $T_2$ are the existence times given by
	Theorem~\ref{thm:vxste}.
	Suppose that
	\[
	u,v\in C_t^\infty(V\times[0,T])
	\cap C([0,T];L^\infty(V))
	\]
	are two solutions of the flow \eqref{Eq:VCLBF} associated with
	$d_0$, satisfying
	\[
	u_i(0)=v_i(0)=0,\qquad i\in V.
	\]
	Then
	\[
	u_i(t)=v_i(t),
	\qquad i\in V,\quad t\in[0,T].
	\]
\end{theorem}

\subsection{Organization of the paper}

The paper is organized as follows.
In Section~2, we recall the angle variation formulas for Euclidean
vertex scalings, derive uniform estimates, and prove short-time
existence of the Euclidean combinatorial Calabi flow.
Section~3 establishes the corresponding results for hyperbolic
vertex scalings.
In Section~4, we use these estimates to prove short-time uniqueness
in both background geometries by Gr\"onwall's inequality.
\\
\\
\textbf{Acknowledgments:} 
The research of Xu Xu is supported by the National Natural Science Foundation of China under Grant No. 12471057 and Hubei Provincial Natural Science Fund for Outstanding Research
Groups under Grant No. 2026AFA008. 
\\
\\
\noindent\textbf{AI Statement:} All mathematical statements and proofs in this paper were independently developed by the authors, who take full responsibility for the content.
AI tools were only used for English-language editing.

	\section{Combinatorial Calabi flow for Euclidean vertex scalings }
	
\subsection{Basic properties of    Euclidean vertex scalings }

\begin{lemma}\label{Lem: VP}\cite{Luo1}
	Let $\triangle ijk$ be a Euclidean triangle with a PL metric $d_0$, and let $u$ be a discrete conformal factor for the triangle.
	Let $\theta^{jk}_{i},\theta^{ik}_{j},\theta^{ij}_{k}$ denote the interior angles at \(i,j,k\), respectively.
	 Then
	\begin{equation*}
		\frac{\partial \theta^{jk}_i }{\partial u_j} =\frac{\partial \theta^{ik}_j }{\partial u_i}  =\frac{1}{2}\cot \theta^{ij}_k .
	\end{equation*}
	\begin{equation*}
	\frac{\partial \theta^{jk}_i }{\partial u_i}  =-\frac{1}{2}(\cot \theta^{ik}_j +\cot \theta^{ij}_k ).
\end{equation*}
\end{lemma}

\begin{remark}\label{remark:KUVP}
	By the definition of combinatorial curvature, if $j\sim i$, then
	\begin{align}\label{Eq: evdkdu1}
		\frac{\partial K_i}{\partial u_j}=\frac{\partial K_j}{\partial u_i}=- \left(\frac{\partial \theta^{jk_1}_i }{\partial u_j}+\frac{\partial \theta^{jk_2}_i
		 }{\partial u_j} \right)
		=-\frac{1}{2}\big(\cot \theta^{ij}_{k_1}+\cot \theta^{ij}_{k_2}\big).
	\end{align}
If $j=i$, then
	\begin{align}\label{Eq: evdkdu2}
		\frac{\partial K_i}{\partial u_i}=&\sum_{j\sim i} \left(\frac{\partial \theta^{jk_1}_i }{\partial u_j}+\frac{\partial \theta^{jk_2}_i }{\partial u_j}\right)
		=\frac{1}{2}\sum_{j\sim i}\big(\cot \theta^{ij}_{k_1}+\cot \theta^{ij}_{k_2}\big),
	\end{align}
	where $\triangle ijk_1$ and $\triangle ijk_2$ are adjacent triangles sharing the edge $\{ij\}$.
	
	Using the variational formulas for combinatorial curvature \eqref{Eq: evdkdu1} and \eqref{Eq: evdkdu2}, we can write the Euclidean combinatorial Calabi flow \eqref{Eq:VCLBF} as
	\begin{equation}\label{Eq:EVCLBF}
		\frac{du_i}{dt}
		=-\sum_{j\in V} \frac{\partial K_i }{\partial u_j}K_j
		=-\sum_{j\sim i} \frac{\partial K_i }{\partial u_j}\left( K_j-K_i \right)
		=\frac{1}{2} \sum_{j\sim i} \big ( \cot\theta_{k_1}^{ij} +\cot\theta_{k_2}^{ij}\big )\left( K_j-K_i \right).
	\end{equation}
\end{remark}

\begin{proposition}\label{Prop: TimeEvolution}
	Let $(S,\mathcal{T})$ be an infinitely triangulated surface with a PL metric $d_0$, and let $u:V\to\mathbb{R}$ be a function such that $d=u*d_0$ is a PL metric. Along the Euclidean combinatorial Calabi flow, the following evolution equations hold.
	
	For each triangle $\triangle ijk\in F$, the interior angles satisfy
	\begin{equation}\label{Eq: AngleEv1}
		\frac{d\theta_i^{jk}}{dt}
		=
		-\frac12
		\big(
		\cot\theta_j^{ik}
		+
		\cot\theta_k^{ij}
		\big)
		\frac{du_i}{dt}
		+
		\frac12\cot\theta_k^{ij}\frac{du_j}{dt}
		+
		\frac12\cot\theta_j^{ik}\frac{du_k}{dt}.
	\end{equation}
	Moreover, the curvatures satisfy
	\begin{equation}\label{Eq: evdkdt1}
		\frac{dK_i}{dt}
		=
		-\sum_{j\in V}\sum_{k\in V}
		\frac{\partial K_i}{\partial u_j}
		\frac{\partial K_j}{\partial u_k}
		K_k.
	\end{equation}
	Finally, for each vertex $i\in V$,
\begin{equation}\label{Eq: SecondDer}
	\begin{aligned}
		\frac{d^2 u_i}{dt^2} 
		&= \frac12 \sum_{j \sim i} \Bigg[
		- \left( \csc^2 \theta_{k_1}^{ij} \frac{d\theta_{k_1}^{ij}}{dt} + \csc^2 \theta_{k_2}^{ij} \frac{d\theta_{k_2}^{ij}}{dt} \right) (K_j - K_i) \\
		&\qquad + \left( \cot \theta_{k_1}^{ij} + \cot \theta_{k_2}^{ij} \right) \left( \frac{dK_j}{dt} - \frac{dK_i}{dt} \right) \Bigg],
	\end{aligned}
\end{equation}
	where $k_1$ and $k_2$ are the two vertices opposite to the edge $\{ij\}$ in the two adjacent triangles.
\end{proposition}

\proof
	The angle evolution equation follows directly from the chain rule and the variation formulas for the Euclidean vertex scalings:
	
	$$
	\frac{\partial\theta_i^{jk}}{\partial u_i}
	=
	-\frac12
	\left(
	\cot\theta_j^{ik}
	+
	\cot\theta_k^{ij}
	\right),
	\qquad
	\frac{\partial\theta_i^{jk}}{\partial u_j}
	=
	\frac12\cot\theta_k^{ij},
	\qquad
	\frac{\partial\theta_i^{jk}}{\partial u_k}
	=
	\frac12\cot\theta_j^{ik}.
	$$
	Thus, \eqref{Eq: AngleEv1} follows.
	For the curvature evolution, by the chain rule and the Euclidean combinatorial Calabi flow
	$$
	\frac{du_j}{dt}
	=
	-\sum_{k\in V}
	\frac{\partial K_j}{\partial u_k}K_k,
	$$
	we obtain
	$$
	\frac{dK_i}{dt}
	=
	\sum_{j\in V}
	\frac{\partial K_i}{\partial u_j}
	\frac{du_j}{dt}
	=
	-\sum_{j\in V}\sum_{k\in V}
	\frac{\partial K_i}{\partial u_j}
	\frac{\partial K_j}{\partial u_k}
	K_k,
	$$
	which gives \eqref{Eq: evdkdt1}.
	Finally, using the explicit form of the Euclidean combinatorial Calabi flow,
	
	$$
	\frac{du_i}{dt}
	=
	\frac12
	\sum_{j\sim i}
	\left(
	\cot\theta_{k_1}^{ij}
	+
	\cot\theta_{k_2}^{ij}
	\right)
	(K_j-K_i),
	$$
	and differentiating with respect to $t$, we obtain
	$$
	\begin{aligned}
		\frac{d^2u_i}{dt^2}
		&=
		\frac12\sum_{j\sim i}
		\Bigg[
		\frac{d}{dt}
		\left(
		\cot\theta_{k_1}^{ij}
		+
		\cot\theta_{k_2}^{ij}
		\right)
		(K_j-K_i)
		\\
		&\qquad
		+
		\left(
		\cot\theta_{k_1}^{ij}
		+
		\cot\theta_{k_2}^{ij}
		\right)
		\left(
		\frac{dK_j}{dt}
		-
		\frac{dK_i}{dt}
		\right)
		\Bigg].
	\end{aligned}
	$$
	Since
	$$
	\frac{d}{dt}\cot\theta_k^{ij}
	=
	-\csc^2\theta_k^{ij}
	\frac{d\theta_k^{ij}}{dt},
	$$
	equation \eqref{Eq: SecondDer} follows.
\qed

\subsection{Short-time existence for  Euclidean vertex scalings }

We first establish several uniform estimates that will be used in the proof of the short-time existence theorem.

\begin{lemma}\label{Lem: UniformFirstDer}
	Suppose that the vertex degrees of $\mathcal T$ are
	bounded by $D$.
	Along the Euclidean combinatorial Calabi flow,
	at every time $t$ for which the PL metric
	$d(t)=u(t)*d_0$ is $\epsilon$-uniformly nondegenerate,
	we have
	\[
	\left|\frac{du_i}{dt}\right|
	\leq \pi D^2\cot\epsilon,
	\qquad i\in V.
	\]
\end{lemma}

\proof
Since the angles of every Euclidean triangle sum to $\pi$ and each angle is bounded below by $\epsilon$, we have
\(
\epsilon\le \theta_i^{jk}\le \pi-2\epsilon.
\)
Hence 
\(
|\cot\theta_i^{jk}| \leq\cot\epsilon.
\)
It follows from the curvature variation formulas
\eqref{Eq: evdkdu1} that, for $j\sim i$,
\[
\left|\frac{\partial K_i}{\partial u_j}\right|	=	\left|\frac{1}{2}\big(\cot \theta^{ij}_{k_1}+\cot \theta^{ij}_{k_2}\big)\right|\leq\cot\epsilon.
\]
Moreover, by the definition of the combinatorial curvature,
$
K_i=2\pi-\sum_{\triangle ijk\in F}\theta_i^{jk}.
$
Since $\deg(i)\leq D$ and $0<\theta_i^{jk}<\pi$, we obtain
$$2\pi - D\pi \le K_i \le 2\pi.$$
Hence, for any vertices $i,j\in V$,
\begin{equation}\label{Eq: kicha}
	|K_j - K_i| \le \pi D.
\end{equation}
Using the explicit form of the Euclidean combinatorial Calabi flow \eqref{Eq:EVCLBF}, we obtain

$$
\begin{aligned}
	\left|\frac{du_i}{dt}\right|
	&\leq
	\sum_{j\sim i}
	\left|\frac{\partial K_i}{\partial u_j}\right||K_j-K_i|
	\leq
	D\cot\epsilon\cdot\pi D
	=\pi D^2\cot\epsilon.
\end{aligned}
$$

\qed

\begin{proposition}\label{Prop: UniformHigherDer}
	Suppose that the vertex degrees of $\mathcal T$ are
	bounded by $D$.
	Along the Euclidean combinatorial Calabi flow,
	at every time $t$ for which the PL metric
	$d(t)=u(t)*d_0$ is $\epsilon$-uniformly nondegenerate,
	we have
	\[
	\left|\frac{dK_i}{dt}\right|
	\leq C_1(\epsilon,D),
	\qquad
	\left|\frac{d^2u_i}{dt^2}\right|
	\leq C_2(\epsilon,D),
	\qquad i\in V,
	\]
	where
	\[
	C_1(\epsilon,D)
	:=4D^2(D+2)\pi\cot^2\epsilon
	\]
	and
	\[
	C_2(\epsilon,D)
	:=2\pi^2D^4(1+\cot^2\epsilon)\cot^2\epsilon
	+2D C_1(\epsilon,D)\cot\epsilon.
	\]
\end{proposition}

\proof

Since the metric is $\epsilon$-uniformly nondegenerate,
$|\cot\theta^{ij}_{k}|\leq\cot\epsilon$ for every triangle.
Thus, \eqref{Eq: evdkdu1} and \eqref{Eq: evdkdu2} give
\[
\left|\frac{\partial K_i}{\partial u_j}\right|
\leq\cot\epsilon \quad (j\sim i),
\qquad
\left|\frac{\partial K_i}{\partial u_i}\right|
\leq D\cot\epsilon.
\]
As the remaining partial derivatives vanish, we obtain
\begin{equation}\label{eq:eukiuj}
	\sum_{j\in V}
	\left|\frac{\partial K_i}{\partial u_j}\right|
	=
	\left|\frac{\partial K_i}{\partial u_i}\right|
	+\sum_{j\sim i}
	\left|\frac{\partial K_i}{\partial u_j}\right|
	\leq 2D\cot\epsilon.
\end{equation}
Moreover, by the definition of the combinatorial curvature,
\[
|K_i|\leq(D+2)\pi.
\]
Thus, using \eqref{Eq: evdkdt1},
\begin{equation}\label{eq:eukidt}
\begin{aligned}
	\left|\frac{dK_i}{dt}\right|
	&\leq
	\sum_{j\in V}
	\left|\frac{\partial K_i}{\partial u_j}\right|
	\sum_{k\in V}
	\left|\frac{\partial K_j}{\partial u_k}\right|
	|K_k|\leq
	4D^2(D+2)\pi\cot^2\epsilon=C_1(\epsilon,D).
\end{aligned}
\end{equation}
We next estimate the second time derivative of $u_i$. 
By the angle evolution formula \eqref{Eq: AngleEv1} in
Proposition~\ref{Prop: TimeEvolution} and
Lemma~\ref{Lem: UniformFirstDer},

\begin{equation}\label{eq:euangledt}
	\begin{aligned}
		\left|\frac{d\theta_i^{jk}}{dt}\right| &\leq \frac12 \left( |\cot\theta_j^{ik}|+|\cot\theta_k^{ij}| \right) \left|\frac{du_i}{dt}\right|+ \frac12|\cot\theta_k^{ij}| \left|\frac{du_j}{dt}\right| + \frac12|\cot\theta_j^{ik}| \left|\frac{du_k}{dt}\right|\\ &\leq 2\pi D^2\cot^2\epsilon. 
	\end{aligned} 	
\end{equation}
Let $\triangle ijk_1$ and $\triangle ijk_2$ be the two triangles adjacent to the edge $\{ij\}$.
Applying Proposition \ref{Prop: TimeEvolution} and taking absolute values, we obtain
\[
\begin{aligned}
	\left|\frac{d^2u_i}{dt^2}\right|
	\leq\frac12\sum_{j\sim i}\Big[
	&\left(
	\csc^2\theta_{k_1}^{ij}
	\left|\frac{d\theta_{k_1}^{ij}}{dt}\right|
	+
	\csc^2\theta_{k_2}^{ij}
	\left|\frac{d\theta_{k_2}^{ij}}{dt}\right|
	\right)
	|K_j-K_i|\\
	&+
	\left(
	|\cot\theta_{k_1}^{ij}|
	+
	|\cot\theta_{k_2}^{ij}|
	\right)
	\left(
	\left|\frac{dK_j}{dt}\right|
	+
	\left|\frac{dK_i}{dt}\right|
	\right)
	\Big].
\end{aligned}
\]
Combining \eqref{eq:euangledt} and the bound
\eqref{eq:eukidt} with
\[
\csc^2\theta\leq 1+\cot^2\epsilon,\qquad
|\cot\theta|\leq\cot\epsilon,\qquad
|K_j-K_i|\leq\pi D,
\]
and $\deg(i)\leq D$, we obtain
\[
\left|\frac{d^2u_i}{dt^2}\right|
\leq
2\pi^2D^4(1+\cot^2\epsilon)\cot^2\epsilon
+2D C_1(\epsilon,D)\cot\epsilon
=C_2(\epsilon,D).
\]
\qed

We next recall the angle estimate for Euclidean vertex scalings.

\begin{lemma}\label{lem: evzl}\cite{JI}
	Let $\triangle ijk$ be a Euclidean triangle with interior angles $\theta_i,\theta_j,\theta_k$ and edge lengths $d_{ij},d_{ik},d_{jk}$. Suppose each angle satisfies $\theta_i,\theta_j,\theta_k \ge \epsilon>0$. Then there exists $\delta=\delta(\epsilon)>0$ such that,
	whenever $\|u\|_{L^\infty}\leq\delta$, the scaled edge lengths
	satisfy the strict triangle inequalities, and the
	corresponding angles satisfy
	\[
	|\widehat\theta_r-\theta_r|\leq\frac{\epsilon}{2},
	\qquad r\in\{i,j,k\}.
	\]
\end{lemma}

\begin{remark}\label{rem: evzl}
	Under the assumptions of Lemma \ref{lem: evzl}, the new PL metric $\tilde d=u*d$ is $\epsilon/2$-uniformly nondegenerate. Indeed,
	$$
	\hat\theta_r\geq
	\theta_r-\left|\hat\theta_r-\theta_r\right|
	\geq\frac{\epsilon}{2},
	\qquad r\in\{i,j,k\}.
	$$
\end{remark}

We next construct a finite exhaustion of the triangulation.
Fix $i_0\in V$, and let $d(i,j)$ denote the combinatorial
distance between $i$ and $j$, namely the minimum number of
edges in a path joining them. Set
\[
V_n=B_n(i_0):=\{j\in V:d(i_0,j)\leq n\},
\qquad n\in\mathbb N.
\]
Since $\mathcal T$ is locally finite and $S$ is connected,
each $V_n$ is finite,
\[
V_n\subset V_{n+1},
\qquad
\bigcup_{n=1}^{\infty}V_n=V.
\]
Let $\mathcal T_n$ consist of all simplices of $\mathcal T$
whose vertices belong to $V_n$, and denote its edge and face
sets by $E_n$ and $F_n$, respectively.
We call $i\in V_n$ an interior vertex if every face incident
to $i$ belongs to $F_n$, and write
\[
\partial V_n=V_n\setminus\operatorname{int}(V_n).
\]
For each fixed vertex $i\in V$, all triangles incident to $i$
belong to $F_n$ for all sufficiently large $n$.
Hence $i\in\operatorname{int}(V_n)$ for all such $n$.

For each $n$, we extend $u^{[n]}$ by zero to $V\setminus V_n$
and compute the curvatures with respect to the metric
$u^{[n]}*d_0$ on the original triangulation $\mathcal T$.
We write $K_i^{[n]}=K_i(u^{[n]})$ and consider the system
\[
\begin{cases}
	\displaystyle
	\frac{du_i^{[n]}}{dt}
	=-\sum_{j\sim i}
	\frac{\partial K_i}{\partial u_j}(u^{[n]})
	\bigl(K_j^{[n]}-K_i^{[n]}\bigr),
	& i\in\operatorname{int}(V_n),\ t>0,\\[6pt]
	u_i^{[n]}(t)=0,
	& i\in\partial V_n,\ t\geq0,\\[4pt]
	u_i^{[n]}(0)=0,
	& i\in V_n.
\end{cases}
\]
Here the curvature derivatives are those of the original
vertex scaling, evaluated at $u^{[n]}$.
The unknowns are the conformal factors at the finitely many
interior vertices. The right-hand side is smooth whenever
the deformed triangles are nondegenerate.
The Picard--Lindel\"of theorem therefore gives a unique local
smooth solution for each $n$.
Since $u^{[n]}$ is identically zero on
$V\setminus\operatorname{int}(V_n)$,
all its time derivatives vanish there.

\textbf{Proof of Theorem~\ref{thm:vxste}(a).}
By Lemma~\ref{lem: evzl}, there exists
$\delta_0=\delta_0(\epsilon)>0$ such that
$u*d_0$ is an $\epsilon/2$-uniformly nondegenerate PL metric
whenever $\|u\|_{L^\infty}\leq\delta_0$.
Set
\[
T_1:=\frac{\delta_0}{2\pi D^2\cot(\epsilon/2)}.
\]
For each $n$, let $\tau_n$ be the maximal time for which
the solution $u^{[n]}$ is defined on $[0,\tau_n)$ such that
$u^{[n]}(t)*d_0$ is a PL metric and
$\|u^{[n]}(t)\|_{L^\infty}<\delta_0$
for all $t\in[0,\tau_n)$.
By Lemma~\ref{lem: evzl}, the metric $u^{[n]}(t)*d_0$
is $\epsilon/2$-uniformly nondegenerate throughout this
interval.
We estimate $\tau_n$ as follows.
For $t\in[0,\tau_n)$ and $i\in V$, using
Lemma~\ref{Lem: UniformFirstDer}, with $\epsilon$ replaced
by $\epsilon/2$, and $u^{[n]}(0)=0$, we have
\[
\begin{aligned}
	|u_i^{[n]}(t)|
	&\leq \int_0^t
	\left|\frac{du_i^{[n]}}{ds}(s)\right|\,ds\leq \pi D^2\cot\frac{\epsilon}{2}\,t.
\end{aligned}
\]
Now if $\tau_n\leq T_1$, then
\[
\limsup_{t\uparrow\tau_n}
\|u^{[n]}(t)\|_{L^\infty}
\leq
\pi D^2\cot\frac{\epsilon}{2}\,\tau_n
\leq \frac{\delta_0}{2}<\delta_0.
\]
By Lemma~\ref{lem: evzl} and the standard continuation
theorem for finite-dimensional ordinary differential
equations, we can extend $u^{[n]}$ to
$[0,\tau_n+\eta_n)$ for some $\eta_n>0$.
By continuity, after decreasing $\eta_n$ if necessary,
the extended solution still satisfies
$\|u^{[n]}(t)\|_{L^\infty}<\delta_0$
throughout this interval.
This contradicts the maximality of $\tau_n$.
Consequently, $\tau_n>T_1$ for every $n$.
Thus all finite-dimensional solutions are defined on
$[0,T_1]$ and satisfy
\[
\|u^{[n]}(t)\|_{L^\infty}
\leq \pi D^2\cot\frac{\epsilon}{2}\,t
\leq \frac{\delta_0}{2},
\qquad 0\leq t\leq T_1.
\]
Their associated metrics are $\epsilon/2$-uniformly
nondegenerate throughout this interval.
To obtain the higher derivative estimates for the finite systems,
	we first apply the chain rule:
	\[
	\left|\frac{dK_i^{[n]}}{dt}\right|
	=
	\left|
	\sum_{j\in\operatorname{int}(V_n)}
	\frac{\partial K_i}{\partial u_j}(u^{[n]})
	\frac{du_j^{[n]}}{dt}
	\right|
	\leq C_1\left(\frac{\epsilon}{2},D\right),
	\qquad i\in V.
	\]
	Here we have used \eqref{eq:eukiuj} and the uniform first
	derivative bound, with $\epsilon$ replaced by $\epsilon/2$.
	The same angle derivative estimates and differentiation
	of the interior equations as in the proof of
	Proposition~\ref{Prop: UniformHigherDer} then give
	\[
	\left|\frac{d^2u_i^{[n]}}{dt^2}\right|
	\leq C_2\left(\frac{\epsilon}{2},D\right),
	\qquad i\in V.
	\]
	Consequently,
	\[
	\|u_i^{[n]}\|_{C^2([0,T_1])}
	\leq\frac{\delta_0}{2}
	+\pi D^2\cot\frac{\epsilon}{2}
	+C_2\left(\frac{\epsilon}{2},D\right),
	\qquad i\in V.
	\]
	This bound is independent of $n$ and $i$.
	
	For each fixed $i\in V$, the sequence $\{u_i^{[n]}\}$
	is therefore uniformly bounded in $C^1([0,T_1])$,
	and its first derivatives are equicontinuous.
	Since $V$ is countable, the Arzel\`a--Ascoli theorem and
	a diagonal argument yield a subsequence $\{n_l\}$ such that
	\[
	u_i^{[n_l]}\longrightarrow u_i^*
	\qquad\text{in }C^1([0,T_1])
	\]
	for every $i\in V$.
	Passing to the limit in the bound for $|u_i^{[n_l]}(t)|$
	and then taking the supremum over $i$, we obtain
	\[
	\|u^*(t)\|_{L^\infty}\leq\frac{\delta_0}{2},
	\qquad 0\leq t\leq T_1.
	\]
	It follows from Lemma~\ref{lem: evzl} that $u^*(t)*d_0$
	is $\epsilon/2$-uniformly nondegenerate on $[0,T_1]$.
	Fix $i\in V$.
	For all sufficiently large $n$, every triangle incident to
	$i$ or to a vertex adjacent to $i$ belongs to $F_n$.
	All curvatures and coefficients in the equation at $i$
	are smooth functions of the conformal factors at finitely
	many vertices of combinatorial distance at most two from $i$.
	The above $C^1$ convergence therefore allows us to pass
	to the limit in this equation and obtain
	
	\[
	\begin{cases}
		\displaystyle
		\frac{du_i^*}{dt}
		=-\sum_{j\sim i}
		\frac{\partial K_i}{\partial u_j}(u^*)
		\bigl(K_j(u^*)-K_i(u^*)\bigr),\\[6pt]
		u_i^*(0)=0,
	\end{cases}
	\]
	for every $i\in V$ and $t\in[0,T_1]$.
	For each fixed $i$, the right-hand side is a smooth
	function of finitely many conformal factors.
	Since all these factors are $C^1$ in time, the equation
	implies that $u_i^*$ is $C^2$ in time.
	Repeating this argument for all vertices gives
	$u_i^*\in C^\infty([0,T_1])$ for every $i\in V$.
\qed

\begin{remark}\label{rem:euclidean-existence-continuity}
	The solution $u^*$ constructed above is continuous in time
	with respect to the $L^\infty$ norm.
	Indeed, the uniform first derivative estimate passes to
	the limit under the above $C^1$ convergence and gives
	\[
	\left|\frac{du_i^*}{dt}\right|
	\leq\pi D^2\cot\frac{\epsilon}{2},
	\qquad i\in V,\quad t\in[0,T_1].
	\]
	Integrating between $s$ and $t$ and taking the supremum
	over $i\in V$, we obtain
	\[
	\|u^*(t)-u^*(s)\|_{L^\infty}
	\leq\pi D^2\cot\frac{\epsilon}{2}\,|t-s|,
	\qquad s,t\in[0,T_1].
	\]
	Together with $u^*(0)=0$, this shows that
	$u^*\in C([0,T_1];L^\infty(V))$.
	Thus the solution obtained in Theorem~\ref{thm:vxste}(a)
	satisfies the continuity assumption in
	Theorem~\ref{thm:VCLBF_uniqueness}.
\end{remark}

\section{Combinatorial Calabi flow for  hyperbolic vertex scalings }

\subsection{Basic properties of  the  hyperbolic vertex scalings }

For each edge $\{ij\}\in E$, let $\triangle ijk_1$ and
$\triangle ijk_2$ denote the two faces incident to $\{ij\}$. For
$r=1,2$, we set
\begin{equation*}
	\Theta_{ij}^{(r)}
	:=
	\frac{
		\theta_i^{jk_r}
		+\theta_j^{ik_r}
		-\theta_{k_r}^{ij}
	}{2}.
\end{equation*}
We also introduce the notation
\begin{equation*}
	B_{ij}
	:=
	\frac{\partial\theta_i^{jk_1}}{\partial u_j}
	+
	\frac{\partial\theta_i^{jk_2}}{\partial u_j},
\end{equation*}
and
\begin{equation*}
	A_i
	:=
	\sum_{\triangle ijk\in F}
	\frac{\partial\operatorname{Area}(\triangle ijk)}
	{\partial u_i}.
\end{equation*}

\begin{lemma}\label{Lem: HVP}
	\cite{GY,XZ,XZ2}
	Let $\triangle ijk$ be a hyperbolic triangle equipped with a
	PH metric $d_0$, and let $u$ denote its discrete conformal factor.
	Write
	$
	d_{st}:=(u*d_0)_{st}
	$
	for the length of the edge $\{st\}$. Let
	$\theta_i^{jk},\theta_j^{ik},\theta_k^{ij}$ denote the interior
	angles at $i,j,k$, respectively. Then
	\begin{equation}\label{Eq:VP1}
		\frac{\partial\theta_i^{jk}}{\partial u_j}
		=
		\frac{\partial\theta_j^{ik}}{\partial u_i}
		=
		\frac{\cosh d_{jk}+\cosh d_{ik}-\cosh d_{ij}-1}
		{\mathcal A(\cosh d_{ij}+1)}
		=
		\frac{1}{2\cosh^2\frac{d_{ij}}{2}}
		\tan
		\frac{
			\theta_i^{jk}
			+\theta_j^{ik}
			-\theta_k^{ij}
		}{2},
	\end{equation}
	where
$
	\mathcal A
	=
	\sinh d_{ik}\sinh d_{ij}\sin\theta_i^{jk}.
$
\end{lemma}

\begin{remark}\label{remark: hkiujkiui}
	By the definition of combinatorial curvature and Lemma \ref{Lem: HVP}, if \(j\sim i\), we obtain
\begin{equation}\label{Eq:hvdkdu1}
		\begin{aligned}
			\frac{\partial K_i}{\partial u_j}	=
			\frac{\partial K_j}{\partial u_i}
			&=
			-\left(
			\frac{\partial\theta_i^{jk_1}}{\partial u_j}
			+
			\frac{\partial\theta_i^{jk_2}}{\partial u_j}
			\right)=
			-\frac{1}{2\cosh^2\frac{d_{ij}}{2}}
			\left(
			\tan\Theta_{ij}^{(1)}
			+
			\tan\Theta_{ij}^{(2)}
			\right)
			=-B_{ij},
		\end{aligned}
	\end{equation}
	where $\triangle ijk_1$ and $\triangle ijk_2$ are the two faces
	in $F$ incident to the edge $\{ij\}$.
	
	For $j=i$, we have
	\begin{equation}\label{Eq:hvdkdu2}
		\begin{aligned}
			\frac{\partial K_i}{\partial u_i}
			&=
			\sum_{\triangle ijk\in F}
			\frac{\partial\bigl(
				\theta_j^{ik}+\theta_k^{ij}
				+\operatorname{Area}(\triangle ijk)
				\bigr)}{\partial u_i}\\
			&=
			\sum_{j\sim i}
			\left(
			\frac{\partial\theta_i^{jk_1}}{\partial u_j}
			+
			\frac{\partial\theta_i^{jk_2}}{\partial u_j}
			\right)
			+
			\sum_{\triangle ijk\in F}
			\frac{\partial\operatorname{Area}(\triangle ijk)}{\partial u_i}\\
			&=
			\sum_{j\sim i}B_{ij}+A_i.
		\end{aligned}
	\end{equation}
Using \eqref{Eq:hvdkdu1} and \eqref{Eq:hvdkdu2}, we can write
the hyperbolic combinatorial Calabi flow \eqref{Eq:VCLBF} as
	\begin{equation}\label{Eq:hVCLBF}
		\begin{aligned}
			\frac{du_i}{dt}
			&=
			-\sum_{j\in V}
			\frac{\partial K_i}{\partial u_j}K_j=
			-\sum_{j\sim i}
			\frac{\partial K_i}{\partial u_j}K_j
			-\frac{\partial K_i}{\partial u_i}K_i=
			\sum_{j\sim i}B_{ij}(K_j-K_i)-A_iK_i.
		\end{aligned}
	\end{equation}
\end{remark}

\begin{remark}\label{remark: hdarea}
	We shall use the following variational identity for the area
	of a hyperbolic triangle due to Glickenstein and Thomas
	\cite[Proposition 9]{GT}:
	\begin{equation}\label{eq:ddarea}
		\frac{\partial \operatorname{Area}(\triangle ijk)}{\partial u_i}
		=
		\frac{\partial\theta_j^{ik}}{\partial u_i}
		\bigl(\cosh d_{ij}-1\bigr)
		+
		\frac{\partial\theta_k^{ij}}{\partial u_i}
		\bigl(\cosh d_{ik}-1\bigr).
	\end{equation}
	Combining \eqref{eq:ddarea} with Lemma~\ref{Lem: HVP},
	we obtain the following explicit expression for $A_i$,
	which also appeared in our previous work \cite{BX}:
	\begin{equation}\label{Eq:A_i_explicit}
		A_i
		=
		\sum_{j\sim i}
		\tanh^2\frac{d_{ij}}{2}
		\left(
		\tan\Theta_{ij}^{(1)}
		+
		\tan\Theta_{ij}^{(2)}
		\right),
	\end{equation}
	where $\triangle ijk_1$ and $\triangle ijk_2$ are the two
	faces incident to the edge $\{ij\}$.
	
	Moreover, since
	$\operatorname{Area}(\triangle ijk)
	=\pi-\theta_i^{jk}-\theta_j^{ik}-\theta_k^{ij}$,
	the identity \eqref{eq:ddarea} gives
	\begin{equation}\label{eq:theta_ui}
		\frac{\partial\theta_i^{jk}}{\partial u_i}
		=
		-\cosh d_{ij}\,
		\frac{\partial\theta_j^{ik}}{\partial u_i}
		-\cosh d_{ik}\,
		\frac{\partial\theta_k^{ij}}{\partial u_i}.
	\end{equation}
\end{remark}

\begin{proposition}\label{Prop:HSecondDer}
	Let $(S,\mathcal{T})$ be an infinitely triangulated surface with a
	PH metric $d_0$, and let $u:V\to\mathbb{R}$ be a function such that
	$d=u*d_0$ is a PH metric. Along the hyperbolic combinatorial Calabi
	flow, for each vertex $i\in V$, we have
	\begin{equation}\label{Eq:HSecondDer}
		\begin{aligned}
			\frac{d^2u_i}{dt^2}
			={}&
			\sum_{j\sim i}
			\frac{dB_{ij}}{dt}(K_j-K_i)
			+
			\sum_{j\sim i}
			B_{ij}\left(\frac{dK_j}{dt}-\frac{dK_i}{dt}\right)
		-
			\frac{dA_i}{dt}K_i
			-
			A_i\frac{dK_i}{dt}.
		\end{aligned}
	\end{equation}
	Moreover,
\begin{equation}\label{Eq:HBt}
	\begin{aligned}
		\frac{d B_{ij}}{dt}
	&=
		-\frac{1}{2}\operatorname{sech}^2\frac{d_{ij}}{2}
		\tanh^2\frac{d_{ij}}{2}
		\left( \frac{du_i}{dt} + \frac{du_j}{dt} \right)
		\left( \tan\Theta_{ij}^{(1)} + \tan\Theta_{ij}^{(2)} \right)\\
		&\quad + \frac{1}{2\cosh^2\frac{d_{ij}}{2}}
		\sum_{r=1}^2 \sec^2\Theta_{ij}^{(r)} \frac{d\Theta_{ij}^{(r)}}{dt}.
	\end{aligned}
\end{equation}
	Furthermore, the curvature evolution is given by
	\begin{equation}\label{Eq:HKt}
		\frac{dK_i}{dt}
		=
		-\sum_{j\in V}\sum_{k\in V}
		\frac{\partial K_i}{\partial u_j}
		\frac{\partial K_j}{\partial u_k}
		K_k.
	\end{equation}
	Finally, we have
	\begin{equation}\label{Eq:HAt}
	\begin{aligned}
		\frac{dA_i}{dt}
		&=
		\sum_{j\sim i} \Bigg[
		\tanh^2\frac{d_{ij}}{2}\operatorname{sech}^2\frac{d_{ij}}{2}
		\left( \tan\Theta_{ij}^{(1)} + \tan\Theta_{ij}^{(2)} \right)
		\left( \frac{du_i}{dt} + \frac{du_j}{dt} \right)\\
		&\quad + \tanh^2\frac{d_{ij}}{2}
		\sum_{r=1}^2 \sec^2\Theta_{ij}^{(r)} \frac{d\Theta_{ij}^{(r)}}{dt}
		\Bigg].
	\end{aligned}
\end{equation}
	
\end{proposition}

\proof By the explicit form of the hyperbolic combinatorial Calabi flow
\eqref{Eq:hVCLBF}, we have
\[
\begin{aligned}
	\frac{d^2u_i}{dt^2}
	={}&
	\sum_{j\sim i}\frac{dB_{ij}}{dt}(K_j-K_i)
	+\sum_{j\sim i}B_{ij}
	\left(\frac{dK_j}{dt}-\frac{dK_i}{dt}\right)-\frac{dA_i}{dt}K_i-A_i\frac{dK_i}{dt},
\end{aligned}
\]
where we have differentiated \eqref{Eq:hVCLBF} and used the product
rule. This proves \eqref{Eq:HSecondDer}.
We compute the time derivative of $B_{ij}$:
\begin{equation*}
	\begin{aligned}
		\frac{d B_{ij}}{dt}
		&=
		\frac{d}{dt}\left[
		\frac{1}{2\cosh^2\frac{d_{ij}}{2}}
		\left( \tan\Theta_{ij}^{(1)} + \tan\Theta_{ij}^{(2)} \right)
		\right] \\
		&=
		\left( \frac{d}{dt}\frac{1}{2\cosh^2\frac{d_{ij}}{2}} \right)
		\left( \tan\Theta_{ij}^{(1)} + \tan\Theta_{ij}^{(2)} \right)\\
		&\quad + \frac{1}{2\cosh^2\frac{d_{ij}}{2}}
		\left( \frac{d}{dt}\tan\Theta_{ij}^{(1)} 
		+ \frac{d}{dt}\tan\Theta_{ij}^{(2)} \right) \\
		&=
		-\frac{1}{2}\operatorname{sech}^2\frac{d_{ij}}{2}
		\tanh^2\frac{d_{ij}}{2}
		\left( \frac{du_i}{dt} + \frac{du_j}{dt} \right)
		\left( \tan\Theta_{ij}^{(1)} + \tan\Theta_{ij}^{(2)} \right)\\
		&\quad + \frac{1}{2\cosh^2\frac{d_{ij}}{2}}
		\sum_{r=1}^2 \sec^2\Theta_{ij}^{(r)} \frac{d\Theta_{ij}^{(r)}}{dt}.
	\end{aligned}
\end{equation*}
Here we have used the identity
	\[
	\frac{d}{dt}
	\left(
	\cosh^{-2}\frac{d_{ij}}{2}
	\right)
	=
	-2\operatorname{sech}^2\frac{d_{ij}}{2}
	\tanh\frac{d_{ij}}{2}
	\frac{d}{dt}\left(\frac{d_{ij}}{2}\right),
	\]
together with the edge-length evolution formula obtained from hyperbolic vertex scalings,
	\[
	\frac{dd_{ij}}{dt}
	=
\tanh\frac{d_{ij}}{2}
	\left(
	\frac{du_i}{dt}
	+
	\frac{du_j}{dt}
	\right).
	\]
	Moreover, 
	we obtain
	\[
	\frac{d\Theta_{ij}^{(r)}}{dt}
	=
	\frac12
	\left(
	\frac{d\theta_i^{jk_r}}{dt}
	+
	\frac{d\theta_j^{ik_r}}{dt}
	-
	\frac{d\theta_{k_r}^{ij}}{dt}
	\right).
	\]
	The angle derivatives are obtained from the chain rule:
	\[
	\frac{d\theta_i^{jk}}{dt}
	=
	\frac{\partial\theta_i^{jk}}{\partial u_i}\frac{du_i}{dt}
	+
	\frac{\partial\theta_i^{jk}}{\partial u_j}\frac{du_j}{dt}
	+
	\frac{\partial\theta_i^{jk}}{\partial u_k}\frac{du_k}{dt}.
	\]
	Here the partial derivatives are given by
	Lemma~\ref{Lem: HVP} and \eqref{eq:theta_ui}.
	Analogous formulas hold for $\theta_j^{ik}$ and $\theta_k^{ij}$.
	For the curvature, the chain rule and \eqref{Eq:hVCLBF} give
\[
\begin{aligned}
	\frac{dK_i}{dt}
	&=
	\sum_{j\in V}
	\frac{\partial K_i}{\partial u_j}\frac{du_j}{dt}=
	-\sum_{j\in V}\sum_{k\in V}
	\frac{\partial K_i}{\partial u_j}
	\frac{\partial K_j}{\partial u_k}K_k,
\end{aligned}
\]
	which proves \eqref{Eq:HKt}.
	Finally, by Remark \ref{remark: hdarea},
\[
A_i
=
\sum_{j\sim i}
\tanh^2\frac{d_{ij}}2
\left(
\tan\Theta_{ij}^{(1)}
+\tan\Theta_{ij}^{(2)}
\right).
\]
Hence,
\[
\begin{aligned}
	\frac{dA_i}{dt}
	={}&
\sum_{j\sim i}
	\tanh^2\frac{d_{ij}}2
	\operatorname{sech}^2\frac{d_{ij}}2
	\left(
	\tan\Theta_{ij}^{(1)}
	+\tan\Theta_{ij}^{(2)}
	\right)
	\left(
	\frac{du_i}{dt}+\frac{du_j}{dt}
	\right)\\
	&+
	\sum_{j\sim i}
	\tanh^2\frac{d_{ij}}2
	\sum_{r=1}^2
	\sec^2\Theta_{ij}^{(r)}
	\frac{d\Theta_{ij}^{(r)}}{dt},
\end{aligned}
\]
where the first term follows from differentiating
$\tanh^2(d_{ij}/2)$ and using the above formula for $dd_{ij}/dt$,
while the second term follows from the chain rule.
\qed

\subsection{Short-time existence for hyperbolic vertex scalings}
As in the Euclidean case, we first establish the uniform
estimates needed to prove short-time existence.

\begin{lemma}\label{Lem: HUniformFirstDer}
	Suppose that the vertex degrees of $\mathcal T$ are
	bounded by $D$.
	Along the hyperbolic combinatorial Calabi flow,
	at every time $t$ for which the PH metric
	$d(t)=u(t)*d_0$ is $\epsilon$-uniformly nondegenerate,
	we have
	\[
	\left|\frac{du_i}{dt}\right|
	\leq \mathcal C_1(\epsilon,D),
	\qquad i\in V,
	\]
	where
	\[
	\mathcal C_1(\epsilon,D)
	:=\pi D^2\cot\epsilon
	+2D(D+2)\pi\cot\epsilon.
	\]
\end{lemma}
\proof
Since the metric is $\epsilon$-uniformly nondegenerate, every
interior angle satisfies
\(
\epsilon\leq\theta_i^{jk}<\pi.
\)
As the angles of a hyperbolic triangle sum to less than $\pi$, we
have
$$
-\frac{\pi}{2}+\epsilon
<
\Theta_{ij}^{(r)}
<
\frac{\pi}{2}-\epsilon,
\qquad r=1,2.
$$
Hence
$$
\left|\tan\Theta_{ij}^{(r)}\right|
\leq\cot\epsilon.
$$
Combining \eqref{Eq:hvdkdu1}  and
$0<\cosh^{-2}(d_{ij}/2)\leq1$, it follows that
\begin{align*}
	|B_{ij}|=\left|\frac{\partial \theta^{jk_1}_i }{\partial u_j}+\frac{\partial \theta^{jk_2}_i }{\partial u_j}\right|
	=\frac{1}{2\cosh^2\frac{d_{ij}}{2}} \left|\tan\Theta_{ij}^{(1)}
	+\tan\Theta_{ij}^{(2)} \right|< \cot\left(\epsilon\right).
\end{align*}
Similarly, by \eqref{Eq:A_i_explicit} and
$0<\tanh^2(d_{ij}/2)<1$,
\begin{align*}
	|A_i|\le \sum_{j\sim i }\tanh^2\frac{d_{ij}}{2}\left|\tan\Theta_{ij}^{(1)}
	+\tan\Theta_{ij}^{(2)} \right|< 2D\cot\left(\epsilon\right).
\end{align*}
Moreover, we have
\[
|K_i|\leq(D+2)\pi,
\qquad
|K_j-K_i|\leq\pi D.
\]
Therefore, by \eqref{Eq:hVCLBF},
\[
\begin{aligned}
	\left|\frac{du_i}{dt}\right|
	&\leq
	\sum_{j\sim i}|B_{ij}|\,|K_j-K_i|
	+
	|A_i|\,|K_i|\\
	&\leq
	\pi D^2\cot\epsilon
	+
	2D(D+2)\pi\cot\epsilon\\
	&=
	\mathcal C_1(\epsilon,D).
\end{aligned}
\]
\qed

\begin{proposition}\label{Prop: HUniformHigherDer}
	Suppose that the vertex degrees of $\mathcal T$ are
	bounded by $D$.
	Along the hyperbolic combinatorial Calabi flow,
	at every time $t$ for which the PH metric
	$d(t)=u(t)*d_0$ is $\epsilon$-uniformly nondegenerate,
	we have
	\begin{equation}\label{Eq:HUniformHigherDer}
		\left|\frac{dK_i}{dt}\right|
		\leq \mathcal C_K(\epsilon,D),
		\qquad
		\left|\frac{d^2u_i}{dt^2}\right|
		\leq \mathcal C_3(\epsilon,D),
		\qquad i\in V,
	\end{equation}
	where
	\begin{align*}
		\mathcal C_K(\epsilon,D)
		&:=16D^2(D+2)\pi\cot^2\epsilon,\\
		\mathcal C_3(\epsilon,D)
		&:=\pi D^2\mathcal C_B(\epsilon,D)
		+4D\cot\epsilon\cdot\mathcal C_K(\epsilon,D)
		+(D+2)\pi\mathcal C_A(\epsilon,D).
	\end{align*}
\end{proposition}

\proof
For $j\sim i$, the hyperbolic curvature variation formulas
\eqref{Eq:hvdkdu1} and \eqref{Eq:hvdkdu2} give
\[
\left|\frac{\partial K_i}{\partial u_j}\right|=|B_{ij}|
\leq \cot\epsilon,
\]
and
\[
\left|\frac{\partial K_i}{\partial u_i}\right|
\leq \sum_{j\sim i}\left|B_{ij}\right|+	\left|A_i\right|\leq3D\cot\epsilon .
\]
Therefore,
\begin{equation}\label{eq:hukiuj}
	\sum_{j\in V}
	\left|\frac{\partial K_i}{\partial u_j}\right|
	=
	\left|\frac{\partial K_i}{\partial u_i}\right|
	+\sum_{j\sim i}
	\left|\frac{\partial K_i}{\partial u_j}\right|
	\leq 4D\cot\epsilon.
\end{equation}
Using \(\lvert K_k\rvert\leq(D+2)\pi\), \eqref{eq:hukiuj}, and the curvature evolution formula \eqref{Eq:HKt}, we obtain
\begin{equation}\label{eq:HKt-bound}
	\begin{aligned}
		\left|\frac{dK_i}{dt}\right|
	&\leq
	\sum_{j\in V}
	\left|\frac{\partial K_i}{\partial u_j}\right|
	\sum_{k\in V}
	\left|\frac{\partial K_j}{\partial u_k}\right|
	|K_k|\leq
	16D^2(D+2)\pi\cot^2\epsilon=\mathcal C_K(\epsilon,D).
	\end{aligned}
\end{equation}
Next, we estimate the time derivatives of the angles. For a
hyperbolic triangle, the hyperbolic cosine law yields
\[
\cosh d_{ij}
=
\frac{
	\cos\theta_k^{ij}
	+
	\cos\theta_i^{jk}\cos\theta_j^{ik}
}{
	\sin\theta_i^{jk}\sin\theta_j^{ik}
}
\leq
\frac{2}{\sin^2\epsilon}.
\]
By \eqref{eq:theta_ui} and Lemma \ref{Lem: HVP},
\begin{equation}\label{eq:theta_ui_abs}
	\left|\frac{\partial \theta^{jk}_i}{\partial u_i}\right|
	\le \cosh d_{ij}\,\left|\frac{\partial \theta^{ik}_j}{\partial u_i}\right|
	+\cosh d_{ik}\,\left|\frac{\partial \theta^{ij}_k}{\partial u_i}\right|	\leq
	\frac{2}{\sin^2\epsilon} \cot\epsilon,
\end{equation}
while
\[
\left|
\frac{\partial\theta_i^{jk}}{\partial u_j}
\right|
+
\left|
\frac{\partial\theta_i^{jk}}{\partial u_k}
\right|
\leq  \cot\epsilon.
\]
Therefore, by the chain rule,
\begin{equation}\label{Eq:HAngleDer}
	\begin{aligned}
		\left|\frac{d\theta_i^{jk}}{dt}\right|
		&\leq
		\sum_{r\in\{i,j,k\}}
		\left|\frac{\partial\theta_i^{jk}}{\partial u_r}\right|
		\left|\frac{du_r}{dt}\right|\leq
		\left(\frac{2}{\sin^2\epsilon}+1\right)
		\mathcal C_1(\epsilon,D)\cot\epsilon=\mathcal C_2(\epsilon,D).
	\end{aligned}
\end{equation}
Consequently,
\[
\left|
\frac{d\Theta_{ij}^{(r)}}{dt}
\right|
\leq
\frac32\mathcal C_2(\epsilon,D),
\qquad r=1,2.
\]
We now estimate the time derivative of $B_{ij}$.
Using \eqref{Eq:HBt} and taking absolute values, we obtain
\begin{equation}\label{eq:HBt-bound}
	\begin{aligned}
		\left|\frac{d B_{ij}}{dt}\right|
		&\le 
		\frac{1}{2}\operatorname{sech}^2\frac{d_{ij}}{2}
		\tanh^2\frac{d_{ij}}{2}
		\left| \frac{du_i}{dt} + \frac{du_j}{dt} \right|
		\left| \tan\Theta_{ij}^{(1)} + \tan\Theta_{ij}^{(2)} \right|\\
		&+ \frac{1}{2\cosh^2\frac{d_{ij}}{2}}
		\sum_{r=1}^2 \sec^2\Theta_{ij}^{(r)} \left|\frac{d\Theta_{ij}^{(r)}}{dt}\right|\\
		&\leq
		2\mathcal C_1(\epsilon,D)\cot\epsilon
		+
		\frac32(1+\cot^2\epsilon)
		\mathcal C_2(\epsilon,D)\\
		&=	\mathcal C_B(\epsilon,D).
	\end{aligned}
\end{equation}
Similarly, using \eqref{Eq:HAt} and taking absolute values, we obtain
\begin{equation}\label{eq:HAt-bound}
	\begin{aligned}
		\left|\frac{dA_i}{dt}\right|
	\le{}&
	\sum_{j\sim i}
	\tanh^2\frac{d_{ij}}2
	\operatorname{sech}^2\frac{d_{ij}}2
	\left|
	\tan\Theta_{ij}^{(1)}
	+\tan\Theta_{ij}^{(2)}
	\right|
	\left|
	\frac{du_i}{dt}+\frac{du_j}{dt}
	\right|\\
	&+
	\sum_{j\sim i}
	\tanh^2\frac{d_{ij}}2
	\sum_{r=1}^2
	\sec^2\Theta_{ij}^{(r)}
	\left|
	\frac{d\Theta_{ij}^{(r)}}{dt}
	\right|\\
	&\leq D\left[
	4\mathcal C_1(\epsilon,D)\cot\epsilon
	+
	3(1+\cot^2\epsilon)
	\mathcal C_2(\epsilon,D)
	\right]\\
	&=
	\mathcal C_A(\epsilon,D).
	\end{aligned}
\end{equation}
Finally, using \eqref{Eq:HSecondDer} and taking absolute values, we obtain
\[
\begin{aligned}
	\left|\frac{d^2u_i}{dt^2}\right|
	\le{}&
	\sum_{j\sim i}
	\left|\frac{dB_{ij}}{dt}\right|
	|K_j-K_i|
	+
	\sum_{j\sim i}
	|B_{ij}|
	\left(
	\left|\frac{dK_j}{dt}\right|
	+
	\left|\frac{dK_i}{dt}\right|
	\right)+
	\left|\frac{dA_i}{dt}\right|
	|K_i|
	+
	|A_i|
	\left|\frac{dK_i}{dt}\right|.
\end{aligned}
\]
Combining \eqref{eq:HKt-bound}, \eqref{eq:HBt-bound},
and \eqref{eq:HAt-bound} with
\[
|K_j-K_i|\leq\pi D,\qquad
|K_i|\leq(D+2)\pi,\qquad
|B_{ij}|\leq\cot\epsilon,\qquad
|A_i|\leq2D\cot\epsilon,
\]
and the degree bound $\deg(i)\leq D$, we obtain
\[
\begin{aligned}
	\left|\frac{d^2u_i}{dt^2}\right|
	&\leq
	\pi D^2\mathcal C_B(\epsilon,D)
	+
	4D\mathcal C_K(\epsilon,D)\cot\epsilon
	+
	(D+2)\pi\mathcal C_A(\epsilon,D)=
	\mathcal C_3(\epsilon,D).
\end{aligned}
\]
\qed

We next recall the angle estimate for hyperbolic vertex
scaling.

\begin{lemma}\label{lem: hvac}\cite{BX}
	For any hyperbolic triangle $\triangle ijk$,  we denote by $\theta_i$, $\theta_j$, and $\theta_k$ the interior angles of the triangle, and by $d_{ij}$, $d_{ik}$ and $d_{jk}$ the edge lengths. Assume
	that $\theta_i,\theta_j,\theta_k \ge \epsilon>0$. Then there exists $\delta=\delta(\epsilon)>0$ such that,
	whenever $\|u\|_{L^\infty}\leq\delta$, the scaled edge lengths
	satisfy the strict triangle inequalities, and the
	corresponding angles satisfy
	\[
	|\widehat\theta_r-\theta_r|\leq\frac{\epsilon}{2},
	\qquad r\in\{i,j,k\}.
	\]
\end{lemma}

\begin{remark}\label{rem: hvzl}
	As in the Euclidean case, under the assumptions of Lemma
	\ref{lem: hvac}, the new PH metric $\tilde d=u*d$ remains
	$\epsilon/2$-uniformly nondegenerate.
\end{remark}

We use the finite exhaustion $\mathcal T_n$ and the zero
extension of $u^{[n]}$ introduced in the Euclidean case.
All curvatures and coefficients below are computed on
$\mathcal T$ with respect to the PH metric $u^{[n]}*d_0$.
Set
\[
K_i^{[n]}=K_i(u^{[n]}),\qquad
B_{ij}^{[n]}=B_{ij}(u^{[n]}),\qquad
A_i^{[n]}=A_i(u^{[n]}).
\]
We consider the finite-dimensional system
\[
\begin{cases}
	\displaystyle
	\frac{du_i^{[n]}}{dt}
	=\sum_{j\sim i}B_{ij}^{[n]}
	\bigl(K_j^{[n]}-K_i^{[n]}\bigr)
	-A_i^{[n]}K_i^{[n]},
	& i\in\operatorname{int}(V_n),\ t>0,\\[6pt]
	u_i^{[n]}(t)=0,
	& i\in\partial V_n,\ t\geq0,\\[4pt]
	u_i^{[n]}(0)=0,
	& i\in V_n.
\end{cases}
\]
As in the Euclidean case, the Picard--Lindel\"of theorem
gives a unique local smooth solution for each $n$.

\textbf{Proof of Theorem~\ref{thm:vxste}(b).}
By Lemma~\ref{lem: hvac}, there exists
$\delta_0=\delta_0(\epsilon)>0$ such that
$u*d_0$ is an $\epsilon/2$-uniformly nondegenerate PH metric
whenever $\|u\|_{L^\infty}\leq\delta_0$.
Set
\[
T_2:=\frac{\delta_0}{2\mathcal C_1(\epsilon/2,D)}.
\]
For each $n$, define $\tau_n$ as in the proof of
Theorem~\ref{thm:vxste}(a), with PL metrics replaced
by PH metrics.
For $0\leq t<\tau_n$, Lemma~\ref{Lem: HUniformFirstDer},
with $\epsilon$ replaced by $\epsilon/2$,
and the fixed zero values outside $\operatorname{int}(V_n)$
give
\[
\left|\frac{du_i^{[n]}}{dt}\right|
\leq\mathcal C_1\left(\frac{\epsilon}{2},D\right),
\qquad i\in V.
\]
Integrating and using $u^{[n]}(0)=0$, we obtain
\[
\|u^{[n]}(t)\|_{L^\infty}
\leq\mathcal C_1\left(\frac{\epsilon}{2},D\right)t,
\qquad 0\leq t<\tau_n.
\]
The same continuation argument as in the proof of
Theorem~\ref{thm:vxste}(a) shows that $\tau_n>T_2$
for every $n$.
Thus all finite-dimensional solutions are defined on
$[0,T_2]$ and satisfy
\[
\|u^{[n]}(t)\|_{L^\infty}
\leq\mathcal C_1\left(\frac{\epsilon}{2},D\right)t
\leq\frac{\delta_0}{2},
\qquad 0\leq t\leq T_2.
\]
Their associated metrics are $\epsilon/2$-uniformly
nondegenerate throughout this interval.
	To apply the estimates in
	Proposition~\ref{Prop: HUniformHigherDer} to the finite
	systems, we use the chain rule and \eqref{eq:hukiuj} to obtain
	\[
	\left|\frac{dK_i^{[n]}}{dt}\right|
	\leq4D\cot\frac{\epsilon}{2}\,
	\mathcal C_1\left(\frac{\epsilon}{2},D\right)
	\leq\mathcal C_K\left(\frac{\epsilon}{2},D\right).
	\]
	The estimates for the time derivatives of $B_{ij}$ and $A_i$
	in that proposition depend only on the uniform angle and
	first derivative bounds, so they also apply here.
	Differentiating the equation at each interior vertex
	therefore gives
	\[
	\left|\frac{d^2u_i^{[n]}}{dt^2}\right|
	\leq\mathcal C_3\left(\frac{\epsilon}{2},D\right).
	\]
	At the remaining vertices the second derivative is zero.
	Consequently,
	\[
	\|u_i^{[n]}\|_{C^2([0,T_2])}
	\leq\frac{\delta_0}{2}
	+\mathcal C_1\left(\frac{\epsilon}{2},D\right)
	+\mathcal C_3\left(\frac{\epsilon}{2},D\right),
	\qquad i\in V.
	\]
	
	As in the Euclidean case, the Arzel\`a--Ascoli theorem and
	a diagonal argument yield a subsequence such that
	\[
	u_i^{[n_l]}\longrightarrow u_i^*
	\qquad\text{in }C^1([0,T_2])
	\]
	for every $i\in V$.
	The limit satisfies $u^*(0)=0$ and
	$\|u^*(t)\|_{L^\infty}\leq\delta_0/2$.
	Thus $u^*(t)*d_0$ remains
	$\epsilon/2$-uniformly nondegenerate.
	For each fixed vertex, the right-hand side is a smooth
	function of finitely many conformal factors.
	We may therefore pass to the limit in the finite-dimensional
	equations and then differentiate successively in time,
	as in the proof of Theorem~\ref{thm:vxste}(a).
	This gives a solution
	$u^*\in C_t^\infty(V\times[0,T_2])$
	of the hyperbolic combinatorial Calabi flow.
	\qed

\begin{remark}\label{rem:hyperbolic-existence-continuity}
	As in Remark~\ref{rem:euclidean-existence-continuity},
	the uniform first derivative estimate gives
	\[
	\|u^*(t)-u^*(s)\|_{L^\infty}
	\leq\mathcal C_1\left(\frac{\epsilon}{2},D\right)|t-s|,
	\qquad s,t\in[0,T_2].
	\]
	Since $u^*(0)=0$, we have
	$u^*\in C([0,T_2];L^\infty(V))$.
	Thus the solution obtained in Theorem~\ref{thm:vxste}(b)
	satisfies the continuity assumption in
	Theorem~\ref{thm:VCLBF_uniqueness}.
\end{remark}

\section{Short-time uniqueness of the combinatorial Calabi flow}

In this section, we prove short-time uniqueness of the combinatorial
Calabi flow in both background geometries. We first use the estimates
established in Sections 2 and 3 to control the difference of the
right-hand sides, and then apply Gr\"onwall's inequality.
Here and below, $\mathcal C_1$ denotes the constant defined
in Lemma \ref{Lem: HUniformFirstDer}, while $\mathcal C_2$, $\mathcal C_B$, and
$\mathcal C_A$ denote the constants defined in the proof of
Proposition~\ref{Prop: HUniformHigherDer}.

\begin{lemma}[Uniform estimate for the right-hand side]\label{rev:lipschitz}
	Let $(S,\mathcal T)$ be an infinitely triangulated surface with vertex
	degree bounded by $D$, and let $d_0$ be an
	$\epsilon$-uniformly nondegenerate initial PL or PH metric.
	Choose $\delta_0=\delta_0(\epsilon)>0$ so that
	$\|w\|_{L^\infty}\leq\delta_0$ implies that $w*d_0$ is
	$\epsilon/2$-uniformly nondegenerate. Set
	\[
	F_i(w)=-\sum_{j\in V}\frac{\partial K_i}{\partial u_j}(w)K_j(w).
	\]
	Then
	\[
	\|F(u)-F(v)\|_{L^\infty}
	\leq C(\epsilon,D)\|u-v\|_{L^\infty}
	\]
	for all $u,v\in L^\infty(V)$ satisfying
	$\|u\|_{L^\infty},\|v\|_{L^\infty}\leq\delta_0$, where
	\begin{equation}\label{eq:calabi-uniqueness-constant}
		C(\epsilon,D)
		:=
		16D^2\cot^2\frac\epsilon2
		+(D+2)\pi
		\frac{2D\mathcal C_B(\epsilon/2,D)
			+\mathcal C_A(\epsilon/2,D)}
		{\mathcal C_1(\epsilon/2,D)}.
	\end{equation}
\end{lemma}

\proof
	Set $h=u-v$ and $w_s=(1-s)v+su$ for $s\in[0,1]$. Then
	\[
	\|w_s\|_{L^\infty}
	\leq(1-s)\|v\|_{L^\infty}+s\|u\|_{L^\infty}
	\leq\delta_0.
	\]
	Thus $w_s*d_0$ is $\epsilon/2$-uniformly nondegenerate for every $s$.
	In the following estimates, all angles, edge lengths, and partial
	derivatives are evaluated at $w_s$, unless otherwise indicated.
	By the definition of combinatorial curvature and the estimates
	\eqref{eq:eukiuj} and \eqref{eq:hukiuj} established in the proofs of
	Propositions~\ref{Prop: UniformHigherDer}
	and~\ref{Prop: HUniformHigherDer},
	with $\epsilon$ replaced by $\epsilon/2$, we have
	
	\begin{equation}\label{eq:usumku}
		\|K(w_s)\|_{L^\infty}\leq(D+2)\pi,
		\qquad
		\sum_{j\in V}\left|
		\frac{\partial K_i}{\partial u_j}(w_s)\right|
		\leq4D\cot\frac\epsilon2.
	\end{equation}
	For each fixed vertex $i$, the curvature $K_i$ depends only on the
	conformal factors at $i$ and its adjacent vertices. Therefore, the
	fundamental theorem of calculus gives
	\[
	K_i(u)-K_i(v)
	=\int_0^1\sum_{j\in V}
	\frac{\partial K_i}{\partial u_j}(w_s)h_j\,ds.
	\]
	It follows that
	\begin{equation}\label{eq:calabi-curvature-difference}
		\|K(u)-K(v)\|_{L^\infty}
		\leq4D\cot\frac\epsilon2\,\|u-v\|_{L^\infty}.
	\end{equation}
	We next estimate the difference of the curvature derivatives.
	
	In the Euclidean case, we apply the angle variation formulas along
	$w_s$, with $d(w_s)_r/ds=h_r$.
	The calculation leading to \eqref{eq:euangledt} then gives
	\begin{equation}\label{eq:euuangleds}
		\left|\frac{d\theta_i^{jk}(w_s)}{ds}\right|
		=\left|\sum_{r\in\{i,j,k\}}
		\frac{\partial\theta_i^{jk}}{\partial u_r}(w_s)h_r\right|
		\leq2\cot\frac\epsilon2\,\|h\|_{L^\infty}.
	\end{equation}
	For each edge $\{ij\}$, let $k_1$ and $k_2$ be the vertices
	opposite to $\{ij\}$ in the two adjacent triangles.
	The Euclidean curvature variation formulas in Remark \ref{remark:KUVP} give
	\[
	\frac{\partial K_i}{\partial u_j}
	=-\frac12\left(
	\cot\theta_{k_1}^{ij}+\cot\theta_{k_2}^{ij}
	\right),
	\qquad j\sim i,
	\]
	and
	\[
	\frac{\partial K_i}{\partial u_i}
	=-\sum_{j\sim i}\frac{\partial K_i}{\partial u_j}.
	\]
Differentiating the identity for $\partial K_i/\partial u_j$ along $w_s$, we obtain
	\[
	\frac d{ds}\left(
	\frac{\partial K_i}{\partial u_j}(w_s)
	\right)
	=
	\frac12\sum_{r=1}^2
	\csc^2\theta_{k_r}^{ij}
	\frac{d\theta_{k_r}^{ij}}{ds}.
	\]
	Using the preceding angle estimate \eqref{eq:euuangleds} and
	$\csc^2\theta\leq1+\cot^2(\epsilon/2)$, we have
	\[
	\begin{aligned}
		\left|
		\frac d{ds}\left(
		\frac{\partial K_i}{\partial u_j}(w_s)
		\right)
		\right|
		&\leq
		\frac12\sum_{r=1}^2
		\csc^2\theta_{k_r}^{ij}
		\left|\frac{d\theta_{k_r}^{ij}}{ds}\right|\leq
		2\cot\frac\epsilon2
		\left(1+\cot^2\frac\epsilon2\right)
		\|h\|_{L^\infty}.
	\end{aligned}
	\]
Differentiating the identity for $\partial K_i/\partial u_i$ along $w_s$ and taking absolute values, we obtain
	\[
	\left|
	\frac d{ds}\left(
	\frac{\partial K_i}{\partial u_i}(w_s)
	\right)
	\right|
	\leq
	\sum_{j\sim i}\left|
	\frac d{ds}\left(
	\frac{\partial K_i}{\partial u_j}(w_s)
	\right)
	\right|.
	\]
	Since $\partial K_i/\partial u_j=0$ whenever
	$j\ne i$ and $j\not\sim i$, it follows that
	\[
	\begin{aligned}
		\sum_{j\in V}\left|
		\frac d{ds}\left(
		\frac{\partial K_i}{\partial u_j}(w_s)
		\right)
		\right|
		&=
		\left|
		\frac d{ds}\left(
		\frac{\partial K_i}{\partial u_i}(w_s)
		\right)
		\right|
		+\sum_{j\sim i}\left|
		\frac d{ds}\left(
		\frac{\partial K_i}{\partial u_j}(w_s)
		\right)
		\right|\\
		&\leq
		2\sum_{j\sim i}\left|
		\frac d{ds}\left(
		\frac{\partial K_i}{\partial u_j}(w_s)
		\right)
		\right|\\
		&\leq
		4D\cot\frac\epsilon2
		\left(1+\cot^2\frac\epsilon2\right)
		\|h\|_{L^\infty},
	\end{aligned}
	\]
	where we have used $\deg(i)\leq D$.
	
	In the hyperbolic case, since $d(w_s)_r/ds=h_r$, the same chain-rule computation
	used to obtain \eqref{Eq:HAngleDer} in the proof of
	Proposition~\ref{Prop: HUniformHigherDer} gives
	\[
	\begin{aligned}
		\left|\frac{d\theta_i^{jk}(w_s)}{ds}\right|
		&\leq
		\left(\frac{2}{\sin^2(\epsilon/2)}+1\right)
		\cot\frac\epsilon2\,\|h\|_{L^\infty}=
		\frac{\mathcal C_2(\epsilon/2,D)}{\mathcal C_1(\epsilon/2,D)}
		\|h\|_{L^\infty},
	\end{aligned}
	\]
	where $\mathcal C_2$ is the constant in the angle derivative estimate
	in Proposition~\ref{Prop: HUniformHigherDer}. Hence
	\[
	\left|\frac{d\Theta_{ij}^{(r)}}{ds}\right|
	\leq
	\frac{3\mathcal C_2(\epsilon/2,D)}{2\mathcal C_1(\epsilon/2,D)}
	\|h\|_{L^\infty},
	\qquad r=1,2.
	\]
	Moreover, the hyperbolic vertex scaling relation gives
	\[
	\frac{dd_{ij}}{ds}
	=\tanh\frac{d_{ij}}2(h_i+h_j).
	\]
	Applying the same computations for $B_{ij}$ and $A_i$
	as in the proof of Proposition \ref{Prop: HUniformHigherDer}, we obtain
	\begin{equation*}
		\begin{aligned}
			\left|\frac{dB_{ij}(w_s)}{ds}\right|
			&\leq
			\left[
			2\cot\frac\epsilon2
			+\frac32\left(1+\cot^2\frac\epsilon2\right)
			\frac{\mathcal C_2(\epsilon/2,D)}{\mathcal C_1(\epsilon/2,D)}
			\right]\|h\|_{L^\infty}=
			\frac{\mathcal C_B(\epsilon/2,D)}{\mathcal C_1(\epsilon/2,D)}
			\|h\|_{L^\infty},
		\end{aligned}
	\end{equation*}
	and
	\[
	\begin{aligned}
		\left|\frac{dA_i(w_s)}{ds}\right|
		&\leq
		D\left[
		4\cot\frac\epsilon2
		+3\left(1+\cot^2\frac\epsilon2\right)
		\frac{\mathcal C_2(\epsilon/2,D)}{\mathcal C_1(\epsilon/2,D)}
		\right]\|h\|_{L^\infty}=
		\frac{\mathcal C_A(\epsilon/2,D)}{\mathcal C_1(\epsilon/2,D)}
		\|h\|_{L^\infty},
	\end{aligned}
	\]
	where the last equalities follow from the definitions of
	$\mathcal C_B$ and $\mathcal C_A$.
	By the hyperbolic curvature variation formulas in Remark \ref{remark: hkiujkiui}, we have
	\[
	\frac{\partial K_i}{\partial u_j}=-B_{ij},
	\qquad j\sim i,
	\]
	and
	\[
	\frac{\partial K_i}{\partial u_i}
	=\sum_{j\sim i}B_{ij}+A_i.
	\]
	Differentiating these identities along $w_s$ gives
	\[
	\frac d{ds}\left(
	\frac{\partial K_i}{\partial u_j}(w_s)
	\right)
	=-\frac{dB_{ij}(w_s)}{ds},
	\qquad j\sim i,
	\]
	and
	\[
	\frac d{ds}\left(
	\frac{\partial K_i}{\partial u_i}(w_s)
	\right)
	=
	\sum_{j\sim i}\frac{dB_{ij}(w_s)}{ds}
	+\frac{dA_i(w_s)}{ds}.
	\]
	Since $\partial K_i/\partial u_j=0$ whenever
	$j\ne i$ and $j\not\sim i$, we obtain
	\begin{equation}\label{eq:uddkiu}
		\begin{aligned}
			\sum_{j\in V}\left|
			\frac d{ds}\left(
			\frac{\partial K_i}{\partial u_j}(w_s)
			\right)
			\right|
			&=
			\left|
			\sum_{j\sim i}\frac{dB_{ij}(w_s)}{ds}
			+\frac{dA_i(w_s)}{ds}
			\right|
			+\sum_{j\sim i}
			\left|\frac{dB_{ij}(w_s)}{ds}\right|\\
			&\leq
			2\sum_{j\sim i}
			\left|\frac{dB_{ij}(w_s)}{ds}\right|
			+\left|\frac{dA_i(w_s)}{ds}\right|\\
			&\leq
			\frac{2D \mathcal C_B(\epsilon/2,D)+\mathcal C_A(\epsilon/2,D)}
			{\mathcal C_1(\epsilon/2,D)}
			\|h\|_{L^\infty},
		\end{aligned}
	\end{equation}
	where the last inequality follows from  $\deg(i)\leq D$.
	
	By the definitions of $\mathcal C_B$ and $\mathcal C_A$, we have
	\[
	\mathcal C_A(\epsilon/2,D)=2D\mathcal C_B(\epsilon/2,D).
	\]
	Consequently,
	\[
	\begin{aligned}
		\frac{2D \mathcal C_B(\epsilon/2,D)+\mathcal C_A(\epsilon/2,D)}
		{\mathcal C_1(\epsilon/2,D)}
		&=
		4D\frac{\mathcal C_B(\epsilon/2,D)}{\mathcal C_1(\epsilon/2,D)}\\
		&=
		4D\cot\frac\epsilon2
		\left[
		2+\frac32\left(1+\cot^2\frac\epsilon2\right)
		\left(\frac{2}{\sin^2(\epsilon/2)}+1\right)
		\right]\\
		&\geq
		4D\cot\frac\epsilon2
		\left(1+\cot^2\frac\epsilon2\right).
	\end{aligned}
	\]
	Therefore, the preceding bound \eqref{eq:uddkiu} for the derivatives of
	$\partial K_i/\partial u_j$ holds in both background geometries.
	Integrating with respect to $s$, we obtain
	\begin{equation}\label{eq:calabi-coefficient-difference}
		\begin{aligned}
			\sum_{j\in V}\left|
			\frac{\partial K_i}{\partial u_j}(u)
			-\frac{\partial K_i}{\partial u_j}(v)\right|
			&\leq\int_0^1\sum_{j\in V}\left|
			\frac d{ds}\left(\frac{\partial K_i}{\partial u_j}(w_s)\right)
			\right|\,ds\\
			&\leq
			\frac{2D \mathcal C_B(\epsilon/2,D)+\mathcal C_A(\epsilon/2,D)}
			{\mathcal C_1(\epsilon/2,D)}\,\|u-v\|_{L^\infty}.
		\end{aligned}
	\end{equation}
	All sums above contain only finitely many nonzero terms for each fixed
	$i$, since $\partial K_i/\partial u_j=0$ whenever $j\ne i$ and
	$j\not\sim i$.
	By the definition of $F_i$, we have
	\[
	\begin{aligned}
		F_i(u)-F_i(v)
		={}&-\sum_{j\in V}\left(
		\frac{\partial K_i}{\partial u_j}(u)
		-\frac{\partial K_i}{\partial u_j}(v)\right)K_j(u)-\sum_{j\in V}\frac{\partial K_i}{\partial u_j}(v)
		\bigl(K_j(u)-K_j(v)\bigr).
	\end{aligned}
	\]
	Combining the preceding uniform bounds \eqref{eq:usumku} with
	\eqref{eq:calabi-curvature-difference} and
	\eqref{eq:calabi-coefficient-difference}, we obtain
	\[
	\begin{aligned}
		|F_i(u)-F_i(v)|
		&\leq\|K(u)\|_{L^\infty}
		\sum_{j\in V}\left|
		\frac{\partial K_i}{\partial u_j}(u)
		-\frac{\partial K_i}{\partial u_j}(v)\right|\\
		&\quad+\sum_{j\in V}\left|
		\frac{\partial K_i}{\partial u_j}(v)\right|
		\|K(u)-K(v)\|_{L^\infty}\\
		&\leq\left[
		(D+2)\pi
		\frac{2D \mathcal C_B(\epsilon/2,D)+\mathcal C_A(\epsilon/2,D)}
		{\mathcal C_1(\epsilon/2,D)}
		+16D^2\cot^2\frac\epsilon2
		\right]\|u-v\|_{L^\infty}.
	\end{aligned}
	\]
	Taking the supremum over $i\in V$ completes the proof.
	\qed

\textbf{Proof of  Theorem~\ref{thm:VCLBF_uniqueness}.}
	Let $\delta_0$ be the constant used in the proof of
	Theorem~\ref{thm:vxste} for the corresponding background geometry,
	and use the same $\delta_0$ in Lemma~\ref{rev:lipschitz}.
	Recall that the existence times were chosen as
	\[
	T_1=\frac{\delta_0}{2\pi D^2\cot(\epsilon/2)}
	\quad\text{in the Euclidean case},
	\]
	and
	\[
	T_2=\frac{\delta_0}{2\mathcal C_1(\epsilon/2,D)}
	\quad\text{in the hyperbolic case}.
	\]
	Let $u$ and $v$ be as in the theorem.
	We first show that
	\[
	\|u(t)\|_{L^\infty},\ \|v(t)\|_{L^\infty}
	\leq\frac{\delta_0}{2},
	\qquad 0\leq t\leq T.
	\]
	Consider $u$.
	Since $u\in C([0,T];L^\infty(V))$ and $u(0)=0$,
	we have $\|u(t)\|_{L^\infty}<\delta_0$ for all sufficiently small $t$.
	Suppose that $t_0\in(0,T]$ is the first time at which
	$\|u(t_0)\|_{L^\infty}=\delta_0$.
	For $0\leq t<t_0$, the metric $u(t)*d_0$ is
	$\epsilon/2$-uniformly nondegenerate.
	Integrating the estimates in
	Lemmas~\ref{Lem: UniformFirstDer} and
	\ref{Lem: HUniformFirstDer}, respectively, and using $u_i(0)=0$,
	we obtain
	\[
	|u_i(t)|
	\leq
	\begin{cases}
		\displaystyle \pi D^2\cot\frac\epsilon2\,t,
		&\text{in the Euclidean case},\\[2mm]
		\mathcal C_1(\epsilon/2,D)t,
		&\text{in the hyperbolic case}.
	\end{cases}
	\]
	Since $T\leq T_1$ or $T\leq T_2$, according to the background geometry,
	it follows that
	\[
	\|u(t)\|_{L^\infty}\leq\frac{\delta_0}{2},
	\qquad 0\leq t<t_0.
	\]
	Continuity at $t_0$ contradicts
	$\|u(t_0)\|_{L^\infty}=\delta_0$.
	Consequently, $\|u(t)\|_{L^\infty}<\delta_0$ throughout $[0,T]$.
	The same first derivative estimates now apply on the whole interval
	and give $\|u(t)\|_{L^\infty}\leq\delta_0/2$.
	The same argument applies to $v$.
	
	The flow equations and the common initial value imply that
	\[
	u_i(t)-v_i(t)
	=\int_0^t\bigl(F_i(u(s))-F_i(v(s))\bigr)\,ds.
	\]
	Applying Lemma~\ref{rev:lipschitz} and taking the supremum over $i\in V$,
	we obtain
	\[
	\|u(t)-v(t)\|_{L^\infty}
	\leq C(\epsilon,D)
	\int_0^t\|u(s)-v(s)\|_{L^\infty}\,ds,
	\]
	where $C(\epsilon,D)$ is the constant defined in
	\eqref{eq:calabi-uniqueness-constant}.
	The function $t\mapsto\|u(t)-v(t)\|_{L^\infty}$ is continuous.
	Gr\"onwall's inequality therefore gives
	\[
	\|u(t)-v(t)\|_{L^\infty}=0,
	\qquad 0\leq t\leq T.
	\]
\qed

\begin{remark}
	The continuity in the $L^\infty$ norm required in
	Theorem~\ref{thm:VCLBF_uniqueness} follows from a uniform
	positive lower bound for the interior angles along the flow.
	Indeed, suppose that $u(0)=0$ and that $u(t)*d_0$ is
	$\epsilon/2$-uniformly nondegenerate for every $t\in[0,T]$.
	Integrating the first derivative estimates in
	Lemmas~\ref{Lem: UniformFirstDer}
	and~\ref{Lem: HUniformFirstDer}
	between $s,t\in[0,T]$ and taking the supremum over $i\in V$,
	we obtain, for the Euclidean and hyperbolic flows respectively,
	\[
	\|u(t)-u(s)\|_{L^\infty}
	\leq \pi D^2\cot\frac{\epsilon}{2}\,|t-s|
	\]
	and
	\[
	\|u(t)-u(s)\|_{L^\infty}
	\leq \mathcal C_1\left(\frac{\epsilon}{2},D\right)|t-s|.
	\]
	Thus $u\in C([0,T];L^\infty(V))$.
\end{remark}

\end{document}